\documentclass{sn-jnl}
\usepackage[dvips]{epsfig}
\usepackage{amsmath,amssymb,euscript,graphicx,amsfonts,verbatim}
\usepackage{enumerate,color}
\usepackage{graphicx}
\usepackage{xcolor}
\usepackage{soul}

\usepackage{hyperref}
\usepackage{enumerate}

\newtheorem{theorem}{Theorem}[section]
\newtheorem{proposition}[theorem]{Proposition}

\newtheorem{corollary}[theorem]{Corollary}

\newtheorem{remark}[theorem]{Remark}
\begin{document}

\title[Parameters of Quotient-polynomial Graphs]{Parameters of Quotient-polynomial Graphs}

\author*[  1]{\fnm{Allen} \sur{Herman}\footnote{The first author's work has been supported by NSERC.}}\email{Allen.Herman@uregina.ca}

\author[   2]{\fnm{Roghayeh} \sur{Maleki}\footnote{The second author's research is supported in part by the Ministry of Education, Science and Sport of Republic of Slovenia (University of Primorska Developmental funding pillar).}}\email{roghayeh.maleki@famnit.upr.si}

\affil[1]{\orgdiv{Department of Mathematics and Statistics}, \orgname{University of Regina},\orgaddress{\street{3737 Wascana Parkway}, \city{Regina}, \state{SK}, \postcode{S4S 0A2}, \country{Canada}}}

\affil[2]{\orgdiv{FAMNIT}, \orgname{University of Primorska}, \orgaddress{ \street{Glagolja\v ska 8},\city{6000 Koper}, \country{Slovenia}}}

\abstract{
Fiol has characterized quotient-polynomial graphs as precisely the connected graphs whose adjacency matrix generates the adjacency algebra of a symmetric association scheme.   We show that a collection of non-negative integer parameters of size $d + \frac{d(d-1)}{2}$ is adequate for describing symmetric association schemes of class $d$ that are generated by the adjacency matrix of their first non-trivial relation.  We use this to generate a database of the corresponding quotient-polynomial graphs that have small valency and up to 6 classes, and among these find new feasible parameter sets for symmetric association schemes with noncyclotomic eigenvalues. }



\maketitle

\section{Introduction} 

Let $(X,S)$ be a finite association scheme of order $n$ with adjacency matrices $\mathcal{A} =\{A_0=I, A_1, \dots, A_d\}$ and rational adjacency algebra $\mathbb{Q}S = span_{\mathbb{Q}}(\mathcal{A})$. 
If $A$ is an $n \times n$ matrix, let $\mathbb{Q}[A]$ denote the subalgebra of $M_n(\mathbb{Q})$ consisting of polynomials in the matrix $A$.  This paper is concerned with the case where $(X,S)$ is a symmetric association scheme and $\mathbb{Q}S = \mathbb{Q}[A_1]$.   Graphs whose adjacency matrix generates a symmetric association scheme have been called {\it quotient-polynomial} graphs (QPGs) by Fiol \cite{Fiol2016}.  When the adjacency matrix of the QPG is one of the $A_j$'s, we will call it a {\it relational} QPG, and assume the relations have been permuted so that this adjacency matrix is $A_1$.   Distance-regular and cyclotomic graphs are both special cases of such relational QPGs, but there are many others (see \cite{FiolPenjic2021}).   

Our original motivation for this project was to find a way to characterize the parameters of relational QPGs in a manner that would be analogous to the intersection array of a distance-regular graph, and then to use the characterization of these parameters to search for feasible cases of parameters for interesting QPGs, in particular ones that could provide examples of symmetric association schemes with noncyclotomic eigenvalues.   Both of these goals were achieved.  We found a meaningful system of parameters for relational quotient-polynomial graphs (which produces their intersection algebras easily), and used it to generate a database of these parameters subject to some feasibility conditions that can be checked efficiently on a computer for ranks up to 7, order up to 250, and bounded valency depending on the rank.  Among these parameter sets we found several examples of feasible parameter sets for possible QPGs with noncyclotomic eigenvalues.   

\vskip 3mm

A referee of the original version of the article remarked that singly-generated symmetric association schemes, including those that result from Fiol's quotient-polynomial graphs, have been separately discussed in recent articles of Ito \cite{Ito24} and of Xia, Tian, Liang, and Koolen \cite{XTLK23}, and suggested we revise this work to take this into account.  To do so, we fix some notation.  Let $(X,S)$ be an association scheme with relations $s_0 = id_X, s_1, \dots, s_d$ and adjacency matrices $A_0=I, A_1, \dots, A_d$.  Ito defines a {\it polynomial association scheme} to be one whose Bose-Mesner algebra $\mathbb{C}S$ is equal to the ring of polynomials $\mathbb{C}[A]$, where $A$ is the adjacency matrix of some graph; i.e.~$A$ is a sum of distinct adjacency matrices $A_i$ in the symmetric association scheme.  We will refer to association schemes satisfying Ito's definition as being {\it polynomial in a graph}.  (Note that the graphs generating symmetric polynomial association schemes are precisely the quotient-polynomial graphs defined by Fiol in \cite{Fiol2016}.)  The same property has been considered (again separately) in a recent paper of Penji\'c and Monzillo \cite{PenjicMonzilloArxiv}, who consider both symmetric and non-symmetric association schemes that are polynomial in a graph.  They extend Fiol's characterization of quotient-polynomial graphs as being those that produce a ``quotient-polynomial'' partition of the vertices and a ``quotient-polynomial'' refinement of the distance matrices of the graph.  The same refinement of distance matrices of the generating graph is also described by Ito in \cite{Ito24}.  Ito also considers the dual property, which we will refer to as being {\it copolynomial with respect to an idempotent}; i.e. if $\mathbb{C}S^{\circ}$ means the Bose-Mesner algebra with respect to the Kronecker product, then $(X,S)$ is copolynomial in an idempotent $E$ if $\mathbb{C}S^{\circ} = \mathbb{C}[E]^{\circ}$, where $E$ is an idempotent in $\mathbb{C}S$ satisfying $EJ=0$.  Xia, Tian, Liang, and Koolen consider symmetric association schemes for which $\mathbb{C}S = \mathbb{C}[A_j]$, where $A_j$ is the adjacency matrix of one of the nontrivial defining relations of the scheme.  They also consider the dual property, about whether the pointwise multiplication algebras $\mathbb{C}S^{o}$ and $\mathbb{C}[E_i]^{o}$ are equal, where $E_i$ is one of the primitive idempotents of the scheme.  We will refer to these properties as being {\it polynomial in some $A_j$} and {\it copolynomial in some $E_i$}, respectively.  The association schemes whose Bose-Mesner algebra is generated by the adjacency matrix $A_1$ of a relational quotient-polynomial graph are precisely those that are polynomial in $A_1$.   The main result of \cite{XTLK23} answers Ito's question about whether $Q$-polynomial association schemes are always polynomial in the affirmative, by proving $Q$-polynomial association schemes, which are by definition polynomial in $E_1$, are always polynomial in $A_1$, where $A_1$ is the adjacency matrix corresponding to the second-largest dual eigenvalue of $E_1$.  They also prove the dual notion, that $P$-polynomial association scheme that is polynomial in $A_1$ by definition, will always be copolynomial in $E_1$, where $E_1$ is the primitive idempotent corresponding to the second-largest eigenvalue of $A_1$.  Xia et. al. show that the polynomial in some $A_j$ and copolynomial in some $E_i$ properties are equivalent for symmetric association schemes of rank up to $4$, and give examples to show that for rank $5$ and higher these two properties are independent.  Penji\'c and Monzillo also give a characterization of rank $4$ association schemes, both symmetric and not, that are polynomial in a graph \cite{PenjicMonzilloArxiv}.  

These new articles ask several questions and give some conjectures concerning relationships between the various notions of polynomial for association schemes.  We have found our database to be consistent with these conjectures. 

\section{A parameter set system}\label{Sec2}

It is well-known that any commutative association scheme has four equivalent parameter sets: (1) the intersection matrices (whose entries are the intersection numbers), (2) the first eigenmatrix $P$, (3) the second eigenmatrix $Q$, and (4) the dual intersection matrices (whose entries are the Krein parameters).    For symmetric association schemes, we will see that there is an easy way to describe a small set of parameters from which all four of equivalent parameter sets can be determined, but it is not usually minimal.  Finding the minimal number of parameters required to describe all coherent configurations of a given rank is an interesting open problem \cite[Problem 8.1]{Problems}, and the symmetric association scheme case is a special case of it where there is a decent upper bound.  For association schemes that are generated by a distance-regular graph, the intersection array gives a minimal subset of parameters that determines the intersection matrices.   In this section, inspired by what happens for distance-regular graphs, we will describe a small set of parameters that determines the set of intersection matrices for an association scheme arising from a relational quotient-polynomial graph; i.e. a symmetric association scheme that is polynomial in $A_1$.  Conversely, every relational quotient-polynomial graph will determine an equivalence class of these small parameter sets.  

For a distance-regular graph, the first intersection matrix $B_1$ corresponding to the adjacency matrix $A_1$ of the graph is tridiagonal, with row sums equal to its valency $k_1$.  The entries of $B_1$ along the main diagonal are denoted by $a_i$, the entries just above the main diagonal are the $b_i$'s and the entries just below the diagonal are the $c_i$'s.  Since the sum of any row of $B_1$ is $k_1$, we have the identity $c_i + a_i + b_i =k_1$ for $i=0,\dots,d$.  This is why, for distance-regular graphs, to recover the parameters for the association scheme it generates, one only requires the intersection array: $[b_0,b_1,\dots,b_{d-1}; c_1, c_2, \dots, c_d]$.  For relational quotient-polynomial graphs, we seek something along the same lines; a small subset of integer parameters from which the intersection matrices can be recovered.  

For symmetric association schemes, there is a straightforward way to obtain a set of parameters that will suffice.  This makes use of two properties of the intersection matrices.  First, it is well-known that every row sum of the intersection matrix $B_j$ for $j=0,1,\dots,d$ is equal to the corresponding valency $k_j$.    The second is the general identity that the intersection numbers of a symmetric association scheme will satisfy 
\begin{align*}
\lambda_{ij\ell} k_{\ell} = \lambda_{i\ell j} k_j, \mbox{ for all } 0 \le i,j,\ell \le d.
\end{align*}
Let $T_s = \frac{s(s+1)}{2}$ be the $s$-th triangular number.

\bigskip
\begin{proposition}\label{paramSAS} 
	All the intersection numbers of a symmetric association scheme of class $d$ can be recovered from 
	\begin{enumerate}[(i)]
		\item the list of valencies $k_1, \dots, k_d$; 
		\item  the $T_{d-1} + T_{d-2} + \dots + T_1$ intersection numbers $\lambda_{i,j,\ell}$ with $1 \le i \le j < \ell \le d$ that lie below the main diagonals of the non-identity intersection matrices. 
	\end{enumerate}
\end{proposition} 

\smallskip
\begin{proof} 
	The intersection numbers below the main diagonal of the first intersection matrix $B_1$ are given, so using the valencies we can find the intersection numbers $\lambda_{1\ell j} = \lambda_{1 j \ell} k_{\ell}/k_j$ that lie above the main diagonal.  Using the row sum identity we can find the missing entries along the main diagonal.  So the valency $k_1$ and the $T_{d-1}$ intersection numbers $\lambda_{1 j \ell}$ with $1 \le j < \ell \le d$ determine $B_1$. 
	
	Moving from $B_{i-1}$ to $B_i$, we know the entries of the $j$-th columns of $B_i$ for $1 \le j \le i-1 \le d-1$, since $\lambda_{i i' \ell} = \lambda_{i' i \ell}$ for $i' < i$.  Using the same identities as above, we can use the next valency $k_i$ and the next $T_{d-i}$ parameters $\lambda_{ij\ell}$ with $i \le j < \ell \le d$ to determine the remaining entries of $B_i$.  When we reach $B_d$, all but the last column of $B_d$ has already been determined, so if we know $k_d$ we can use the row sum condition to determine the last column of $B_d$.  The Proposition follows.  
\end{proof} 

\medskip
To illustrate Proposition~\ref{paramSAS}, consider the case of a symmetric association scheme with $4$ classes.  Then, Proposition~\ref{paramSAS} tells us the only parameters we need are the ones indicated here: 
$$ B_1 = \begin{bmatrix} 
	0 & k_1 & 0 & 0 & 0 \\ 
	1 & \ast & \ast & \ast & \ast \\ 
	0 & \lambda_{112} & \ast & \ast & \ast \\ 
	0 & \lambda_{113} & \lambda_{123} & \ast & \ast \\ 
	0 & \lambda_{114} & \lambda_{124} & \lambda_{134} & \ast \end{bmatrix},\,\,\,\,\,\,\,
B_2 = \begin{bmatrix}
    0 & 0 & k_2 & 0 & 0 \\ 
	0 & \ast & \ast & \ast & \ast \\ 
	1 & \ast & \ast & \ast & \ast \\ 
	0 & \ast & \lambda_{223} & \ast & \ast \\ 
	0 & \ast & \lambda_{224} & \lambda_{234} & \ast \end{bmatrix},$$\\
$$B_3 = \begin{bmatrix} 
	0 & 0 & 0 & k_3 & 0 \\ 
	0 & \ast & \ast & \ast & \ast \\ 
	0 & \ast & \ast & \ast & \ast \\ 
	1 & \ast & \ast & \ast & \ast \\ 
	0 & \ast & \ast & \lambda_{334} & \ast \end{bmatrix},\,\,\,\,\,\,\,	
B_4 = \begin{bmatrix} 0 & 0 & 0 & 0 & k_4 \\ 
	0 & \ast & \ast & \ast & \ast \\ 
	0 & \ast & \ast & \ast & \ast \\ 
	0 & \ast & \ast & \ast & \ast \\ 
	1 & \ast & \ast & \ast & \ast \end{bmatrix}. $$ 

\medskip	
As is well-known, for distance-regular graphs the intersection matrices are determined entirely by the first intersection matrix $B_1$, and there is a unique natural ordering of $\textbf{B}=\{B_0,B_1,\dots,B_d\}$ for which $B_1$ is tridiagonal.  The non-zero off-diagonal entries of $B_1$ give us the intersection array $[b_0,b_1,\dots,b_{d-1};c_1,c_2,\dots,c_d]$, from which $B_1$ and the orthogonal polynomials $f_i(x)$ with $B_i = f_i(B_1)$ can be determined.   For relational quotient-polynomial graphs, since the association scheme is polynomial in $A_1$, we are interested to know if the list of valencies and the entries below the diagonal of $B_1$ should be enough to determine the association scheme.  With this in mind, we propose the following parameter array pattern for relational quotient-polynomial graphs of class $d$ and valency $k_1$: 
\begin{align}\label{array}
 \quad [[k_1,\dots,k_d],[\lambda_{112},\dots,\lambda_{11d}; \lambda_{123}, \dots, \lambda_{12d}; \lambda_{134}, \dots, \lambda_{1,d-1,d}]].
\end{align}
We will restrict ourselves to only consider parameter arrays of class $d$ and valency $k_1$ that determine a non-negative $(d+1) \times (d+1)$ integral intersection matrix $B_1$ with row sums equal to $k_1$ using the algorithm of Proposition \ref{paramSAS}; we will refer to these as {\it QPG parameter arrays}.   Recall that the set of intersection matrices $\textbf{B}$ of an association scheme is required to be the basis of a standard integral table algebra (SITA) (see \cite{AFM}).    We will say that a QPG parameter array with class $d$ is {\it valid} if the intersection matrix $B_1$ that it determines has nonnegative integer entries and the polynomial ring $\mathbb{Q}[B_1]$ has a $\mathbb{Q}$-basis $\mathbf{B} = \{B_0=I, B_1, B_2, \dots, B_d\}$ that produces structure constants that are compatible with the list of valencies and parameters defined by the array.    We will consider a pair of valid QPG parameter arrays to be {\it equivalent} if the intersection matrices they determine are conjugate by a $(d+1) \times (d+1)$ permutation matrix whose corresponding permutation lies in $Sym(\{2,\dots,d\})$ (i.e. it fixes the first two rows and columns indexed by $0$ and $1$).  

\bigskip
\begin{theorem}\label{QPGarrays}
	There is a one-to-one correspondence between the set of equivalence classes of valid QPG parameter arrays of class $d$ and valency $k_1$ and the exact isomorphism classes of standard integral table algebras of rank $d+1$ that are generated by an element of valency $k_1$.   In particular, every relational quotient-polynomial graph of valency $k_1$ that generates a symmetric association scheme of class $d$ can be associated with a unique equivalence class of valid QPG parameter arrays with class $d$ and valency $k_1$. 
\end{theorem} 

\medskip
\begin{proof} 
Each valid QPG parameter array in (\ref{array}) determines an intersection matrix $B_1$ whose adjacency algebra $\mathbb{Q}[B_1]$ has a standard integral table algebra basis $\textbf{B}$.  This implies that there are polynomials $f_i(x)$ for $i\in \{2,\ldots,d\}$ of degree $\le d+1$ for which the basis element $B_i$ is equal to $f_i(B_1)$.  To find these polynomials, consider the first rows $v_i$ of $B_1^i$ for $i=0,1,\dots, d$.  The first rows of the $B_i$ are equal to $k_i e_i$, where $e_i$ is the elementary standard basis vector with a $1$ in position $i$ for $i\in\{0,1,\dots, d\}$.  If $f_i(x) = \sum_{\ell=0}^d a_{i,\ell} x^{\ell}$, then $B_i = f_i(B_1)$ implies $k_i e_i = \sum_{\ell=0}^d a_{i,\ell} v_{\ell}$ for $i\in \{2,\ldots,d\}$.  If the QPG parameter set is valid, then the coefficients of $f_i(x)$ are found by expressing the $k_i e_i$ in the basis $\{v_0,v_1,\dots,v_d\}$ of $\mathbb{Q}^{d+1}$ coming from the first rows of the powers $B_1^i$ for $i\in\{0,1,\dots,d\}$.  
 
Conversely, if we start with the set of regular matrices of a standard integral table algebra generated by an element $B_1$ of valency $k_1$, then $B_1$ determines a valid QPG parameter set of class $d$, which is unique up to permuting $\{B_2, \dots, B_d\}$.
\end{proof} 

Since the standard integral table algebra above is generated by $B_1$, every $B_i \in \textbf{B}$ must be a rational polynomial in $B_1$.  To find the polynomials $f_i(x)$ for which $B_i = f_i(B_1)$, we can proceed as in the first part of the above proof, and write the vectors $k_ie_i$ as linear combinations of the first row vectors $v_0, v_1, \dots, v_d$ of the powers $B_1^d$.  If $k_ie_i = \sum_j a_{i,j} v_j$, then it must be the case that $B_i = \sum_j a_{i,j} B_1^j$.  

\bigskip
\begin{corollary} 
	A valid QPG parameter array of class $d$ determines a unique set of $d+1$ orthogonal polynomials in a single variable. 
\end{corollary} 

\medskip
A natural question to ask is the following: if a SITA is polynomial in $B_1$, is the SITA determined by the eigenvalues of $B_1$?  The answer to this is no.  We found several pairs of rank $4$ valid QPG parameter arrays for which the resulting $B_1$ has the same minimal polynomial: 

\smallskip
\noindent $[[5,5,1],[2,0;5]]$ and $[[5,15,3],[1,0;5]]$: 

both have $\mu_{B_1}(x) = (x-5)(x+1)(x^2-5)$,

\smallskip
\noindent $[[7,7,1],[3,0;7]]$, $[[7,14,2],[2,0;7]]$, and $[[7,35,5],[1,0;7]$: 

all three have $\mu_{B_1}(x) = (x-7)(x+1)(x^2-7)$.  

\smallskip
\noindent $[[8,6,12],[4,2;2]$ and $[[8,16,2],[3,0;8]]$: 

both have $\mu_1(x) = (x-8)(x-2)(x+1)(x+4)$.  

\smallskip
\noindent The last pair above are both realized as association schemes.  So the eigenvalues of $A_1$ alone do not determine the parameters of an association scheme that is polynomial in $A_1$. The above examples do not tell us if the spectrum of $A_1$ is enough to determine the SITA, because the eigenvalues of the $B_1$'s in these examples will have different multiplicities when considered as eigenvalues of the corresponding $A_1$'s.

\smallskip
For distance-regular graphs, the distance partition imposes a natural ordering of the $A_i$'s, since the degree of $f_i(x)$ is $i$ for $i=1,\dots,d$.  Under this ordering, the pattern of valencies is unimodal: $k_1 \le k_2 \le \dots \le k_t \ge k_{t+1} \ge \dots \ge k_d$ for some $t$, 
the $c_i$'s are increasing: $c_1 \le c_2 \le \dots \le c_d$, 
and the $b_i$'s are decreasing: $k_1 = b_0 \ge b_1 \ge \dots \ge b_{d-1}$.   For relational quotient-polynomial graphs, it is not clear that there will be a preferred way to order the basis elements.  When generating our database, we only made two (safe) assumptions: $\lambda_{112}>0$ and $k_3 \le k_4 \le \dots \le k_d$.  As observed by Fiol in \cite{Fiol2016}, and again by Ito \cite{Ito24} for symmetric association schemes that are polynomial in a graph, the fact that the adjacency algebra is polynomial in $A_1$ imposes a partial order on the set of adjacency matrices.   $A_0$ is at distance $0$, $A_1$ is at distance $1$, $A_i$ is at distance $2$ if $i>1$ and $A_1^2 \circ A_i \ne 0$ (i.e. $i>1$ and $\lambda_{11i}>0$), $A_i$ is at distance $3$ if $A_i$ is not at distance $< 3$ and $A_1^3 \circ A_i \ne 0$, etc.  Following Ito \cite{Ito24}, this defines a distance partition $\Lambda$ of $\{0,1,\dots,d\}$ with $\Lambda_0 = \{0\}$, $\Lambda_1 = \{1\}$, and $\Lambda_k = \{ i : A_i \mbox{ is at distance } k \}$ for $2 \le k \le D$, where $D$ is the diameter of the graph with adjacency matrix $A_1$.  Distance-regular graphs are the extreme case of this, as they have $|\Lambda_i|=1$ for all $i=0,1,\dots,D=d$.   Ito conjectures in \cite{Ito24} that there is an absolute constant $N$ such that whenever a polynomial association scheme has its first $N$ distance partition elements of size $|\Lambda_i| = 1$,  the association scheme will be $P$-polynomial.  

\section{Generating a database of feasible QPG parameter arrays}

Now that we have a convenient system (\ref{array}) for expressing parameters of relational quotient-polynomial graphs, our next goal for this project is to produce a database that gives representatives for the equivalence classes of valid QPG parameter arrays of class $d$ and valency $k_1$ up to a given order $n$ that satisfy the known feasibility conditions for their quotient-polynomial graphs to exist.  We start by generating parameter arrays for the given class $d$ and valency $k_1$ with a given order $n$ that produce a first intersection matrix $B_1$ with non-negative integer entries and all row sums equal to $k_1$. Then we calculate the minimal polynomial $\mu_1(x)$ of $B_1$, and check that it has degree $(d+1)$.  

Next, we finish checking that our QPG parameter array is valid, by solving for a SITA basis of $\mathbb{Q}[B_1]$.  We calculate the powers $B_1^i$ for $i\in \{2,\ldots,d\}$ and let $v_i \in \mathbb{Q}^{d+1}$ be the first row of $B_1^i$ for $i\in\{0,1,\dots,d\}$.  We check that the $\mathbb{Q}$-span of $\{v_0,v_1,\dots,v_d\}$ has dimension $d+1$.  Assuming it does, proceeding as in the proof of Theorem \ref{QPGarrays}, we let $e_i$ be the elementary standard basis of $\mathbb{Q}^{d+1}$ for $i\in\{0,1,\dots,d\}$, and express the $k_i e_i$ as linear combinations in the basis $\{v_0, v_1, \dots, v_d\}$: 
$$ k_i e_i = \sum_{\ell=0}^d a_{i,\ell} v_{\ell}, \mbox{ for } i\in\{0,1,\dots,d\}.$$
Since $k_i e_i$ is naturally the first row of $B_i$ for $i\in\{0,1,\dots,d\}$, and the solution to the $a_{i,\ell}$'s is unique for each $i$, if we set $f_i(x) = \sum_{\ell=0}^d a_{i,\ell} x^{\ell}$, then the matrices $B_i = f_i(B_1)$ for $i=0,1,\dots,d$ will be the elements of the only possible standard table algebra basis $\mathbf{B} = \{B_0,B_1,\dots,B_d\}$ for $\mathbb{Q}[B_1]$ that can be associated with this QPG parameter array.  We need to check directly that these $B_i$ have nonnegative integer entries and that $\mathbf{B} = \{B_0,B_1,\dots,B_d\}$ is a table algebra basis.   

If we do obtain a SITA basis for $\mathbb{Q}[B_1]$, we can quickly check the following feasibility conditions for a symmetric association scheme: 
\begin{enumerate}[(i)]

\item {\bf The handshaking lemma condition:} If any of the $k_i$ that are larger than $1$ are odd, then the entries $(B_i)_{ji}$ for $j>0$ must be even.  
\item \label{fc} {\bf Quick check if $B_1$ allows integral multiplicities:} Find the irreducible factorization $\mu_1(x) = (x-k_1) g_1(x) \cdots g_s(x)$.  Let $g_t(x) = x^{n_t} + g_{t,1} x^{n_t-1} + \dots$ for $t\in\{1,\dots,s\}$, and find non-negative integers $m_1,\dots,m_t$ for which 
	\begin{equation*}
		1+n_1m_1+\dots+n_tm_t = n := 1+k_1+\dots+k_d
	\end{equation*}
and 
\begin{equation}\label{eq}
	k_1 - m_1 g_{1,1} -\dots -m_s g_{s,1}=0.
\end{equation}		
The  equality (\ref{eq}) must be satisfied because it is precisely the column orthogonality relation applied to the $B_1$-column of the first eigenmatrix $P$.   
At least one such list of prospective integral multiplicities needs to exist.  Since we are only checking that one column of $P$ admits integral multiplicities, we also check that these satisfy further requirements necessary for them to be the actual multiplicities of the SITA: the {\it integrality of the Frame number $F$ and its compatibility with the Discriminant of $\mu_1(x)$}.  It is well-known that the {\it Frame number} of a commutative association scheme must be a positive integer.  In \cite{Hanaki2022}, Hanaki has recently observed that $F$ also needs to be a square times the product of the discriminants of the fields occurring in the Wedderburn decomposition of $\mathbb{Q}S$. From feasibility condition~(\ref{fc}), we can calculate our Frame number 
$$ F(S) = n^{d+1} \frac{k_1 \cdots k_d}{m_1^{n_1} \cdots m_t^{n_t}}. $$
Since the degree of the minimal polynomial of $B_1$ is $d+1$, and the powers of $B_1$ lie in $\mathbb{Z}S$, we have that $\mathbb{Z}[B_1]$ is an integral subalgebra of $\mathbb{Z}S$.  Therefore, the discriminant $D$ of the polynomial $\mu_1(x)$ must be equal to $F(S)$ times a perfect square integer.  So, we also check that $D/F(S)$ is a perfect square. 
\end{enumerate}
We repeat the next step for each such list until we find a list of integral multiplicities that passes both steps.  
\begin{enumerate}[(iii)]
\item \label{C1} {\bf Check for integrality of standard trace on powers of $A_1$:} Given a list of potential integral multiplicities found in feasibility condition~(\ref{fc}) , we check for the existence of a {\it standard integral trace}.  Calculate the degree $n$ polynomial 
$$ f(x) = (x-k_1) g_1(x)^{m_1} g_2(x)^{m_2} \dots g_t(x)^{m_t}, $$
and let $C$ be its companion matrix.  Then, calculate $\frac{1}{n} tr(C^i)$ for $i\in\{0,1,\dots,d-1\}$, and check that these values are non-negative integers.  Since $C$ is a conjugate of the (unknown) matrix $A_1$, and the standard feasible trace of the adjacency algebra $\mathbb{Q}S$ is equal to $\frac{1}{n}$ times the usual trace on $n \times n$ matrices, these values are the values of the standard feasible trace on the powers of $A_1^i$, so they must be non-negative integers. 
\end{enumerate}

\medskip
\begin{remark} {\rm In generating our database the SITAs produced by a valid QPG array often have splitting fields that are not easily accommodated in GAP, either because they are noncyclotomic extensions of $\mathbb{Q}$ or they are cyclotomic with a large conductor that cannot be found quickly.  This is why we do not try to produce the eigenmatrices during the sieve.  The quick check that $B_1$ admits integral multiplicities~(\ref{fc}) does not guarantee the other $B_i$ will admit the same integral multiplicities, so it produces many ``false positives".  In the database, we have done what we could to check for integral multiplicities, by computing $P$ and $Q$ when possible to do so, and other times by using the irreducible factorization of the characteristic polynomials of $B_2, \dots, B_d$ and see if their roots satisfy the column orthogonality relations corresponding to a choice of multiplicities compatible with $B_1$.  
}\end{remark}

Our database implementation \cite{HMprogram} produces lists of valid QPG parameter arrays that survive the above sieving process, and indicates if they pass the basic feasibility checks for symmetric association schemes.  These include the handshaking lemma, integral multiplicities, nonnegativity of Krein parameters, and the absolute bound condition (see \cite{HM2023}).  Of course, there are more feasibility conditions for symmetric association schemes that can be checked, but it is a challenge to implement them for ranks larger than $4$.   

As of this writing the database \cite{HMprogram} includes all QPG parameter arrays of order up to $250$ with any valency for rank $4$, valency up to $30$ for rank $5$, valency up to $15$ for rank $6$, and valency up to $12$ for rank $7$.  We have also found all QPG parameter arrays with rank $4$ and $5$ with arbitrary order and valency up to $8$.  There is some overlap of this database with existing databases of distance-regular graphs and $Q$-polynomial association schemes.  We now know from the result in \cite{XTLK23} that each $Q$-polynomial scheme corresponds to a relational quotient polynomial graph.   We have been able to use GAP to confirm that some of feasible QPG arrays of orders between 35 and 60 in our database are realized by association schemes.  It is hoped that this data help to inform future attempts at extending the classification of small association schemes.   

Xia et. al. gave an infinite family of rank $6$ imprimitive examples of symmetric association schemes which are polynomial (in some $A_j$) but not copolynomial (in some $E_i$) \cite[Example 3]{XTLK23}, and indicated Bannai and Zhao had found primitive examples showing these properties are independent of each other.  Our database includes two imprimitive examples, both realizable as association schemes, that are polynomial in $A_1$ but not copolynomial in any of their $E_i$'s: $[[9,2,12,18],[9,0,3;0,0;6]]$ (of order $42$) and $[[8,1,18,24],[8,0,2;0,0;6]]$ (of order $52$).  We note that both of these examples are copolynomial in an idempotent.   The examples given in Xia et al.'s paper, while not copolynomial in any of their $E_j$'s, are copolynomial in their idempotent $E_3 + E_5$, so they are also copolynomial in the sense of Ito's definition.  At this point we do not know of a symmetric association scheme that is polynomial in a graph but not copolynomial in any of its idempotents. 

\section{New parameter sets for symmetric association schemes with noncyclotomic eigenvalues} 

In this last section, we present parameter sets for prospective relational quotient-polynomial graphs with noncyclotomic eigenvalues that emerged from our searches, and indicate which ones pass the short list of feasibility conditions mentioned in the previous section.    

The examples reported here extend the lists of parameter sets of small rank $5$ examples found in \cite{HM2023} and include new examples of rank $6$ and $7$.  The only other such examples we know of are the examples of feasible pseudocyclic parameter sets of large prime order for ranks $5$ and $6$ were found earlier by Hanaki and Teranishi \cite{Hanaki}.  To date, every feasible parameter set for a commutative association scheme with noncyclotomic eigenvalues that has been found can be associated with the existence of a relational quotient-polynomial graph. 

\newpage
\begin{table}[h!] 
\caption{Rank $5$ parameter sets with noncyclotomic eigenvalues, $n \le 250$, $k_1 \le 30$: }
		\begin{tabular}{|c|c|c|}
			\hline
			Order&Parameter Array&Status\footnotemark[1] \\
			\hline
		$35$&$[[4,12,12,6],[1,0,0;1,2;2]]$& nr-xC \cite{HM2023} \\ \hline
		$45$&$[[8,8,24,4],[1,2,0;1,2;6]]$& F \cite{HM2023} \\ \hline
		$76$&$[[12,18,36,9],[2,0,4;4,4;4]]$ & nr-xF \\ \hline
    	$76$ & $[[18,18,36,3],[2,5,0;5,6;12]],$ & F \cite{HM2023} \\ \hline
    	$88$&$[[14,35,35,3],[4,0,14;6,0;0]]$& nr-xF \cite{HM2023} \\ \hline
		$93$&$[[12,30,30,20],[2,2,0;6,3;6]]$& F \cite{HM2023} \\ \hline
		$112$&$[[15,30,36,30],[4,0,2;5,1;6]]$& F \\ \hline
        $116$&$[[19,1,19,76],[19,6,2;0,0;3]]$& nr-xF \cite{HM2023} \\ \hline
		$119$&$[[16,48,48,6],[1,3,0;7,8;8]]$& F \\ \hline
		$120$&$[[17,17,34,51],[2,0,4;3,1;6]]$& F \\ \hline
		$133$&$[[12,24,48,48],[2,1,0;3,2;4]]$& F \\ \hline
		$133$&$[[18,36,72,6],[4,2,0;4,6;12]]$& F \\ \hline	
		$135$&$[[8,56,56,14],[1,0,0;3,4;4]]$& F \\ \hline		
		$190$&$[[18,36,99,36],[1,0,5;4,0;11]]$& F \cite{HM2023} \\ \hline
		$190$&$[[18,54,108,9],[3,1,0;5,6;12]]$& F \\ \hline	
		$209$&$[[10,90,90,18],[1,0,0;4,5;5]]$& F \\ \hline	
		$210$&$[[11,99,66,33],[1,0,0;6,3;6]]$& nr-xF \\ \hline					
		\end{tabular}
\footnotetext{In the third column, the entry reports the status of its feasibility and/or realizability as a QPG, and refers to \cite{HM2023} if this SITA has been reported earlier: \\ nr-xC: not realizable because not in classification; \\ nr-xF: not realizable because it fails one of the easy-to-implement feasibility conditions; \\ F: passes easy-to-implement feasibility conditions. } 
\end{table}

\bigskip
\begin{table}[h!] 
\caption{Rank $6$ parameter sets with noncyclotomic eigenvalues, $n \le 250$, $k_1 \le 16$: }
	\begin{tabular}{|c|c|c|}
		\hline
		Order&Parameter Array&Status\\
		\hline
	       $40$& $[[6,12,1,8,12],[2,6,0,0;0,3,1;0,0;2]]$ & nr-xF \\ \hline 
           $63$ & $[[10,10,2,20,20],[3,0,1,1;0,2,1;0,1;4]]$ & F \\ \hline
           $70$ & $[[8,24,1,12,24],[2,8,0,0;0,4,2;0,0;2]]$ & F \\ \hline
		   $78$ & $[[12,16,1,24,24],[3,12,2,0;0,4,2;0,0;6]$ & F \\ \hline
           $80$ & $[[14,28,1,8,28],[2,14,0,4;0,7,3;0,0;2]]$ & nr-xF \\ \hline
           $84$ & $[[9,27,2,18,27],[2,9,0,0;0,6,3;0,0;0]]$ & nr-xF \\ \hline
           $90$ & $[[16,16,1,8,48],[6,16,0,2;0,4,2;0,0;2]]$ & F \\ \hline
           $105$ & $[[12,18,2,36,36],[2,12,2,0;0,0,2;0,0;7]]$ & F \\ \hline
           $105$ & $[[12,24,2,18,48],[3,12,0,0;0,4,3;0,0;3]]$ & F \\ \hline 
           $105$ & $[[16,16,16,24,32],[4,2,2,1;2,0,4;4,2;3]]$ & F \\ \hline
           $108$ & $[[10,40,1,16,40],[2,10,0,0;0,5,3;0,0;2]]$ & F \\ \hline
           $120$ & $[[9,27,2,27,54],[2,9,0,0;0,5,1;0,0;1]]$ & F \\ \hline 
		   $198$&$[[7,14,1,35,140],[2,7,0,0;0,2,0;0,0;1]]$& nr-xF \\ \hline	
           $208$ & $[[14,84,1,24,84],[2,14,0,0;0,7,5;0,0;2]]$ & F \\ \hline	
           $228$ & $[[9,27,2,63,126],[2,9,0,0;0,3,0;0,0;2]]$ & nr-xF \\ \hline 
	\end{tabular}
\end{table}

\bigskip
\begin{table}[h!] 
\caption{Rank $7$ parameter sets with noncyclotomic eigenvalues, $n \le 250$, $k_1 \le 12$: }
	\begin{tabular}{|c|c|c|}
		\hline
		Order&Parameter Array&Status\\
		\hline
		$36$ & $[[4,8,1,6,8,8],[1,4,0,0,0;0,0,3,0;0,0,0;0,3;1]]$ & nr-xF \\ \hline
        $44$ & $[[12,6,4,6,6,9],[10,0,0,6,4;0,2,0,0;0,2,4;4,4;0]]$ & nr-xF \\ \hline
       $100$& $[[6,24,1,20,24,24],[1,6,0,0,0;0,0,5,0;0,0,0;0,5;1]]$ & nr-xF \\ \hline
       $100$& $[[12,24,1,14,24,24],[5,12,0,0,0;0,0,7,0;0,0,0;0,7;5]]$ & nr-xF \\ \hline
       $126$ & $[[15,20,2,8,40,40],[9,15,0,0,0;0,0,3,0;0,0,0;0,3;12]]$ & nr-xF \\ \hline
       $164$& $[[12,40,1,30,40,40]][3,12,0,0,0;0,4,6,0;0,0,0;0,6;6]]$ & F \\ \hline
       $196$& $[[8,48,1,42,48,48],[1,8,0,0,0;0,0,7,0;0,0,0;0,7;1]]$ & nr-xF \\ \hline
       $220$& $[[12,24,1,14,84,84],[5,12,0,0,0;0,0,2,0;0,0,0;0,2;10]$ & nr-xF \\ \hline
	\end{tabular}
\end{table}

\bigskip
As of this writing, we have made attempts to construct the symmetric association schemes for the feasible parameter sets listed above with orders up to $76$, but so far these attempts have been inconclusive.   

Based on our exhaustive computer searches, we reach the following conclusions for relational quotient polynomial graphs of small rank that have noncyclotomic eigenvalues. 

\medskip
\begin{theorem} 
\begin{enumerate} 
\item All rank $5$ relational quotient polynomial graphs with valency $7$ or less have cyclotomic eigenvalues.  There are two feasible rank $5$ QPG parameter arrays with valency $8$ whose associated association schemes would have noncyclotomic eigenvalues, if they exist. 
\item All rank $6$ relational quotient polynomial graphs with valency up to $7$ and order up to $250$ have cyclotomic eigenvalues.  There are feasible rank $6$ QPG parameter arrays with valency $8$ and $9$ and order less than $250$ whose associated association schemes would have noncyclotomic eigenvalues, if they exist.   
\item All rank $7$ relational quotient polynomial graphs with valency up to $11$ and order up to $250$ have cyclotomic eigenvalues.  There is one feasible rank $7$ QPG parameter arrays with valency $12$ and order less than $250$ whose associated association scheme would have noncyclotomic eigenvalues, if it exists. 
\end{enumerate} 
\end{theorem}  

\bigskip
{\small
{\bf Data availability statement:} Data supporting this article is available at  

\url{https://github.com/RoghayehMaleki/QPGdatabase-}.
}

\end{document}